\documentclass{amsart}
\usepackage{amssymb}
\usepackage{amsmath}
\usepackage{amsfonts}

\newtheorem{theorem}{Theorem}[section]

\newtheorem{corollary}[theorem]{Corollary}
\newtheorem{definition}[theorem]{Definition}
\newtheorem{proposition}[theorem]{Proposition}
\newtheorem{remark}[theorem]{Remark}
\numberwithin{theorem}{section}
\input{tcilatex}

\begin{document}
\title{ON MORITA EQUIVALENCE OF GROUP\ ACTIONS\ ON LOCALLY\ \ $C^{\ast }$%
-ALGEBRAS}
\author{MARIA\ JOI\c{T}A}
\thanks{
This research was partially supported by grant CEEX-code PR-D11-PT00-48/2005
and partially by CNCSIS\ grant A1065 /2006 from the Romanian Ministry of
Education and Research.}
\maketitle

\begin{abstract}
In this paper, we prove that two continuous inverse limit actions $\alpha $
and $\beta $ of a locally compact group $G$ on the locally $C^{\ast }$%
-algebras $A$ and $B$ are strongly Morita equivalent if and only if there is
a locally $C^{\ast }$-algebra $C$ such that $A$ and $B$ appear as two
complementary full corners of $C$ and there is a continuous action $\gamma $
of $G$ on $C$ which leaves $A$ and $B$ invariant such that $\gamma
|_{A}=\alpha $ and $\gamma |_{B}=\beta $. This generalizes a result of
Combes, \textit{Proc. London Math. Soc}. 49(1984), 289-306.
\end{abstract}

\section{Introduction}

Locally $C^{\ast }$-algebras are generalizations of $C^{\ast }$-algebras.
Instead of being by a single $C^{\ast }$-norm, the topology on a locally $%
C^{\ast }$-algebra is defined by a directed family of $C^{\ast }$-seminorms.
Such many concepts as group action on a $C^{\ast }$-algebra, crossed product
of a $C^{\ast }$-algebra by a group action, Hilbert $C^{\ast }$-module,
adjointable module morphism, group action on a Hilbert $C^{\ast }$-module
can be defined in a natural way in the context of locally $C^{\ast }$%
-algebras. The proofs are not always straightforward.

Phillips \cite{11} introduced the notion of action ( inverse limit action)
of a locally compact group $G$ on a metrizable locally $C^{\ast }$-algebra $%
A $ and defined the crossed product of $A$ by a continuous inverse limit
action of $G$ on $A$. In \cite{6}, we proved a version of the Takai duality
theorem for crossed products of locally $C^{\ast }$-algebras by continuous
inverse limit actions. The concept of strong Morita equivalence for locally $%
C^{\ast }$-algebras was introduced in \cite{4}. In \cite{7}, we introduce
the notion of strong Morita equivalence on the set of group actions on
locally $C^{\ast }$-algebras and prove that it is an equivalence relation.
Also, we prove that the crossed products of locally $C^{\ast }$-algebras
associated with two strongly Morita equivalent continuous inverse limit
actions are strongly Morita equivalent.

This paper is organizes as follows. In Section 2 we recall some facts about
Hilbert modules over locally $C^{\ast }$-algebras and ( continuous inverse
limit ) actions of a locally compact group $G$ on a Hilbert module $E$ over
a locally $C^{\ast }$-algebra $A$. In Section 3 we prove that any (
continuous inverse limit ) action of $G$ on a full Hilbert $A$-module $E$
induces a ( continuous inverse limit ) action of $G$ on the linking algebra $%
\mathcal{L}(E)$ of $E$, Proposition 3.1. Also we prove that two continuous
inverse limit actions $\alpha $ and $\beta $ of a locally compact group $G$
on the locally $C^{\ast }$-algebras $A$ and $B$ are strongly Morita
equivalent if and only if there is a locally $C^{\ast }$-algebra $C$ such
that $A$ and $B$ appear as two complementary full corners of $C$ and there
is a continuous inverse limit action $\gamma $ of $G$ on $C$ which leaves $A$
and $B$ invariant such that $\gamma |_{A}=\alpha $ and $\gamma |_{B}=\beta $%
, Theorem 3.5. This generalizes a result of Combes \cite{2} and is a version
of Theorem 3.3 in \cite{4} for group actions on locally $C^{\ast }$-algebras.

\section{Preliminaries}

A locally $C^{\ast }$-algebra is a complete Hausdorff complex topological $%
\ast $-algebra $A$ whose topology is determined by its continuous $C^{\ast }$%
-seminorms in the sense that a net $\{a_{i}\}_{i\in I}$ converges to $0$ in $%
A$ if and only if the net $\{p(a_{i})\}_{i}$ converges to $0$ for all
continuous $C^{\ast }$-seminorm $p$ on $A$. The term of \textquotedblright
locally $C^{\ast }$ -algebra\textquotedblright\ is due to Inoue \cite{3}.

The set $S(A)$ of all continuous $C^{\ast }$-seminorms on $A$ is directed
with the order $p\geq q$ if $p\left( a\right) \geq q\left( a\right) $ for
all $a\in A$. For each $p\in S(A),$ $\ker p=\{a\in A;p(a)=0\}$ is a
two-sided $\ast $-ideal of $A$ and the quotient algebra $A/\ker p$, denoted
by $A_{p}$, is a $C^{\ast }$-algebra in the $C^{\ast }$-norm induced by $p$.
The canonical map from $A$ to $A_{p}$ is denoted by $\pi _{p}$. For $p,q\in
S(A)$ with $p\geq q$ there is a canonical surjective morphism of $C^{\ast }$%
-algebras $\pi _{pq}:A_{p}\rightarrow A_{q}$ such that $\pi _{pq}(\pi
_{p}(a))=\pi _{q}(a)$ for all $a\in A$. Then $\{A_{p};\pi _{pq}\}_{p,q\in
S(A),p\geq q}$ is an inverse system of $C^{\ast }$-algebras and moreover,
the locally $C^{\ast }$-algebras $A$ and $\lim\limits_{\underset{p\in S(A)}{%
\leftarrow }}A_{p}$ are isomorphic.

An approximate unit for $A$ is an increasing net of positive elements $%
\{e_{i}\}_{i\in I}$ in $A$ such that $p(e_{i})\leq 1$ for all $p\in S(A)$
and for all $i\in I$, and $p\left( ae_{i}-a\right) +p(e_{i}a-a)\rightarrow 0$
for all $p\in S(A)$ and for all $a\in A$. Any locally $C^{\ast }$-algebra
has an approximate unit.

A morphism of locally $C^{*}$-algebras is a continuous morphism of $*$
-algebras. Two locally $C^{*}$-algebras $A$ and $B$ are isomorphic if there
is a bijective map $\Phi :A\rightarrow B$ such that $\Phi $ and $\Phi ^{-1}$
are morphisms of locally $C^{*}$-algebras.

Hilbert modules over locally $C^{*}$-algebras are generalizations of Hilbert 
$C^{*}$-modules by allowing the inner-product to take values in a locally $%
C^{*}$-algebra rather than in a $C^{*}$-algebra.

\begin{definition}
A pre -Hilbert$\ A$-module is a complex vector space$\ E$\ which is also a
right $A$-module, compatible with the complex algebra structure, equipped
with an $A$-valued inner product $\left\langle \cdot ,\cdot \right\rangle
:E\times E\rightarrow A\;$which is $\mathbb{C}$ -and $A$-linear in its
second variable and satisfies the following relations:

\begin{enumerate}
\item $\left\langle \xi ,\eta \right\rangle ^{\ast }=\left\langle \eta ,\xi
\right\rangle \;\;$for every $\xi ,\eta \in E;$

\item $\left\langle \xi ,\xi \right\rangle \geq 0\;\;$for every $\xi \in E;$

\item $\left\langle \xi ,\xi \right\rangle \geq 0\;\;$for every $\xi \in
E;\left\langle \xi ,\xi \right\rangle =0\;$\ if and only if $\xi =0.$
\end{enumerate}

We say that $E\;$is a Hilbert $A$-module if $E\;$is complete with respect to
the topology determined by the family of seminorms $\{\overline{p}%
_{E}\}_{p\in S(A)}\;$where $\overline{p}_{E}(\xi )=\sqrt{p\left(
\left\langle \xi ,\xi \right\rangle \right) },\xi \in E$.\smallskip
\end{definition}

Any locally $C^{*}$ -algebra $A$ is a Hilbert $A$ -module in a natural way.

A Hilbert $A$ -module $E$ is full if the linear space $\left\langle
E,E\right\rangle \;$generated by $\{\left\langle \xi ,\eta \right\rangle
,\;\xi ,\eta \in E\}$ is dense in $A$

Let $E\;$be a Hilbert $A$-module.\ For $p\in S(A),\;\ker \overline{p}%
_{E}=\{\xi \in E;\overline{p}_{E}(\xi )=0\}\;$is a closed submodule of $E\;$%
and $E_{p}=E/\ker \overline{p}_{E}\;$is a Hilbert $A_{p}$-module with $(\xi
+\ker \overline{p}_{E}{})\pi _{p}(a)=\xi a+\ker \overline{p}_{E}{}\;$and $%
\left\langle \xi +\ker \overline{p}_{E}{},\eta +\ker \overline{p}%
_{E}{}\right\rangle =\pi _{p}(\left\langle \xi ,\eta \right\rangle ).$\ The
canonical map from $E\;$onto $E_{p}$ is denoted by $\sigma _{p}.$ For $%
p,q\in S(A),\;p\geq q\;$there is a canonical morphism of vector spaces $%
\sigma _{pq}\;$from $E_{p}\;$onto $E_{q}\;$such that $\sigma _{pq}(\sigma
_{p}(\xi ))=\sigma _{q}(\xi ),\;\xi \in E.\;$Then $\{E_{p};A_{p};\sigma
_{pq},\pi _{pq}\}_{p,q\in S(A),p\geq q}$ is an inverse system of Hilbert $%
C^{\ast }$-modules in the following sense: $\sigma _{pq}(\xi
_{p}a_{p})=\sigma _{pq}(\xi _{p})\pi _{pq}(a_{p}),\xi _{p}\in E_{p},a_{p}\in
A_{p};$ $\left\langle \sigma _{pq}(\xi _{p}),\sigma _{pq}(\eta
_{p})\right\rangle =\pi _{pq}(\left\langle \xi _{p},\eta _{p}\right\rangle
),\xi _{p},\eta _{p}\in E_{p};$ $\sigma _{pp}(\xi _{p})=\xi _{p},\;\xi
_{p}\in E_{p}\;$and $\sigma _{qr}\circ \sigma _{pq}=\sigma _{pr}\;$if $p\geq
q\geq r,\ $and $\lim\limits_{\underset{p\in S(A)}{\leftarrow }}E_{p}$ is a
Hilbert $A$-module which can be identified with\ $E$.

The set $L(E)$ of all adjointable $A$-module morphisms from $E$ into $E$ is
a locally $C^{\ast }$-algebra with topology defined by the family of
seminorms $\{\widetilde{p}_{L(E)}\}_{p\in S(A)},$ where $\widetilde{p}%
_{L(E)}(T)=\left\Vert (\pi _{p})_{\ast }(T)\right\Vert _{L(E_{p},)},$ $T\in
L(E)\;$and $(\pi _{p})_{\ast }(T)(\xi +\ker \overline{p}_{E})=T(\xi )+\ker 
\overline{p}_{F},$ $\xi \in E.$ Moreover, $\{L(E_{p});$ $\;(\pi _{pq})_{\ast
}\}_{p,q\in S(A),p\geq q},$ where $(\pi _{pq})_{\ast }:L(E_{p})\rightarrow
L(E_{q})$ is a morphism of $C^{\ast }$-algebras defined by $(\pi
_{pq})_{\ast }(T_{p})(\sigma _{q}(\xi ))=\sigma _{pq}(T_{p}(\sigma _{p}(\xi
)))$, is an inverse system of $C^{\ast }$-algebras, and $\lim\limits_{%
\underset{p\in S(A)}{\leftarrow }}L(E_{p})$ can be identified with $L(E)$.

For $\xi ,\eta \in E$\ \ we consider the rank one homomorphism $\theta
_{\eta ,\xi }$\ from $E$\ into $E$\ defined by $\theta _{\eta ,\xi }(\zeta
)=\eta \left\langle \xi ,\zeta \right\rangle .$\ Clearly, $\theta _{\eta
,\xi }\in L(E)$\ and $\theta _{\eta ,\xi }^{\ast }=\theta _{\xi ,\eta }$.
\smallskip The linear subspace of $L(E)$ spanned by $\{\theta _{\eta ,\xi
};\xi ,\eta \in E\},$ denoted by $\Theta (E)$, is a two sided $\ast $-ideal
of $L(E)$. The closure of $\Theta (E)$ in $L(E)$ is denoted by $K\left(
E\right) $.

Let $E$ be a full Hilbert $A$ -module. Here we recall some facts about the
linking algebra of $E$ from \cite{5}

The direct sum $A\oplus E$ of the Hilbert $A$-modules $A$ and $E$ is a full
Hilbert $A$-module with the action of $A$ on $A\oplus E$ defined by 
\begin{equation*}
\left( A\oplus E,A\right) \ni \left( a\oplus \xi ,b\right) \mapsto \left(
a\oplus \xi \right) b=ab\oplus \xi b\in A\oplus E
\end{equation*}%
and the inner product defined by 
\begin{equation*}
\left( A\oplus E,A\oplus E\right) \ni \left( a\oplus \xi ,b\oplus \eta
\right) \mapsto \left\langle a\oplus \xi ,b\oplus \eta \right\rangle
=a^{\ast }b+\left\langle \xi ,\eta \right\rangle \in A.
\end{equation*}%
Moreover, for each $p\in S(A),$ the Hilbert $A_{p}$-modules $\left( A\oplus
E\right) _{p}$ and $A_{p}\oplus E_{p}$ can be identified \cite{8}. Then the
locally $C^{\ast }$-algebras $L(A\oplus E)$ and $\lim\limits_{\underset{p\in
S(A)}{\leftarrow }}$ $L(A_{p}\oplus E_{p})$ can be identified \cite{10}.

Let $a\in A$, $\xi \in E$, $\eta \in E$ and $T\in K(E)$. The map $L_{a,\xi
,\eta ,T}:A\oplus E\rightarrow A\oplus E$ defined by 
\begin{equation*}
L_{a,\xi ,\eta ,T}(b\oplus \zeta )=(ab+\left\langle \xi ,\zeta \right\rangle
)\oplus (\eta b+T(\zeta ))
\end{equation*}%
is an element in $L(A\oplus E)$. The locally $C^{\ast }$- subalgebra of $%
L(A\oplus E)$ generated by 
\begin{equation*}
\left\{ L_{a,\xi ,\eta ,T};a\in A,\xi \in E,\eta \in E,T\in K(E)\right\}
\end{equation*}%
is denoted by $\mathcal{L}(E)$ and it is called the linking algebra of $E$.

By Lemma III 3.2 in \cite{9}, we have 
\begin{equation*}
\mathcal{L}(E)=\lim\limits_{\underset{p\in S(A)}{\leftarrow }}\overline{%
\left( \pi _{p}\right) _{\ast }(\mathcal{L}(E))},
\end{equation*}%
where $\overline{\left( \pi _{p}\right) _{\ast }(\mathcal{L}(E))}$ means the
closure of the vector space $\left( \pi _{p}\right) _{\ast }(\mathcal{L}(E))$
in $L(A_{p}\oplus E_{p})$. Let $p\in S(A)$. From 
\begin{equation*}
\left( \pi _{p}\right) _{\ast }(L_{a,\xi ,\eta ,T})=L_{\pi _{p}(a),{\sigma
_{p}}(\xi ),{\sigma _{p}}(\eta ),(\pi _{p})_{\ast }(T)}
\end{equation*}%
for all $a,b\in A,$ for all $\xi ,\eta ,\zeta \in E,$ and taking into
account that $\mathcal{L}(E_{p}),$ the linking algebra of $E_{p}$, is
generated by 
\begin{equation*}
\{L_{\pi _{p}(a),{\sigma _{p}}(\xi ),{\sigma _{p}}(\eta ),(\pi _{p})_{\ast
}(T)};a\in A,\xi \in E,\eta \in E,T\in K(E)\}
\end{equation*}%
since $\pi _{p}(A)=A_{p},$ $\sigma _{p}(E)=E_{p},$ and $\overline{\left( \pi
_{p}\right) _{\ast }\left( K(E)\right) }=K(E_{p})$, we conclude that 
\begin{equation*}
\mathcal{L}(E)=\lim\limits_{\underset{p\in S(A)}{\leftarrow }}\mathcal{L}%
(E_{p}).
\end{equation*}%
Moreover, since $\mathcal{L}(E_{p})=K(A_{p}\oplus E_{p})$ \cite{12} and the
locally $C^{\ast }$-algebras $K(A\oplus E)$ and $\lim\limits_{\underset{p\in
S(A)}{\leftarrow }}K(A_{p}\oplus E_{p})$ can be identified, the linking
algebra of $E$ coincides with $K(A\oplus E).$

Here we recall some facts about actions of a locally compact group $G$ on
Hilbert modules from \cite{7}.

Let $A$ and $B$ be two locally $C^{\ast }$-algebras, let $E$ be a Hilbert $A$%
-module and let $F$ be a Hilbert $B$-module.

\begin{definition}
A morphism of Hilbert modules from $E$ to $F$ is a map $u:E\rightarrow F$
with the property that there is a morphism of locally $C^{\ast }$-algebras $%
\alpha :A\rightarrow B$ such that%
\begin{equation*}
\left\langle u\left( \xi \right) ,u\left( \eta \right) \right\rangle =\alpha
\left( \left\langle \xi ,\eta \right\rangle \right)
\end{equation*}%
for all $\xi ,\eta \in E$. An isomorphism of Hilbert modules is a bijective
map $u:E\rightarrow F$ such that $u$ and $u^{-1}$ are morphisms of Hilbert
modules.
\end{definition}

If $u:E\rightarrow F$ is a morphism of Hilbert modules and $\alpha
:A\rightarrow B$ is a morphism of locally $C^{\ast }$-algebras such that $%
\left\langle u\left( \xi \right) ,u\left( \eta \right) \right\rangle =\alpha
\left( \left\langle \xi ,\eta \right\rangle \right) $, then $u$ is a
continuous linear map and $u\left( \xi a\right) =u\left( \xi \right) \alpha
\left( a\right) $ for all $a\in A$ and for all $\xi \in E$. Moreover, if $u$
is an isomorphism of Hilbert modules and the Hilbert modules $E$ and $F$ are
full, then $\alpha $ is an isomorphism of locally $C^{\ast }$-algebras.

For a Hilbert $A$-module $E$, 
\begin{equation*}
\text{Aut}(E)=\{u:E\rightarrow E;u\text{ is an isomorphism of Hilbert
modules }\}
\end{equation*}

is a group.

\begin{definition}
Let $G$ be a locally compact group. An action of $G$ on $E$ is a morphism of
groups $g\mapsto u_{g}$ from $G$ to Aut$(E).$

The action $g\mapsto u_{g}$ of $G$ on $E$ is continuous if the map $G\times
E\ni \left( g,\xi \right) \mapsto u_{g}\left( \xi \right) \in E$ is jointly
continuous.

An action $g\mapsto u_{g}$ of $G$ on $E$ is an inverse limit action if we
can write $E$ as an inverse limit of Hilbert $C^{\ast }$-modules $%
\lim\limits_{\underset{\lambda \in \Lambda }{\leftarrow }}E_{\lambda }$ in
such a way that for each $g\in G,$ $u_{g}=\lim\limits_{\underset{\lambda \in
\Lambda }{\leftarrow }}u_{g}^{\lambda }$, where $g\mapsto u_{g}^{\lambda }$
is an action of $G$ on $E_{\lambda },$ $\lambda \in \Lambda .$
\end{definition}

If $g\mapsto u_{g}$ is an inverse limit action of $G$ on $E$, then $%
E=\lim\limits_{\underset{\lambda \in \Lambda }{\leftarrow }}E_{\lambda }$
and $u_{g}=\lim\limits_{\underset{\lambda \in \Lambda }{\leftarrow }%
}u_{g}^{\lambda }$ for each $g\in G$, where $g\mapsto u_{g}^{\lambda }$ is
an action of $G$ on $E_{\lambda },$ $\lambda \in \Lambda $. Let $\lambda \in
\Lambda $. Since $g\mapsto u_{g}^{\lambda }$ is an action of $G$ on $%
E_{\lambda },$ 
\begin{equation*}
\left\Vert u_{g}^{\lambda }\left( \sigma _{\lambda }(\xi \right) \right\Vert
_{E_{\lambda }}=\left\Vert \sigma _{\lambda }(\xi \right\Vert _{E_{\lambda }}
\end{equation*}%
for each $\xi \in E,$ and for all $g\in G$ \cite{1,2}. This implies that 
\begin{equation*}
\overline{p}_{\lambda }(u_{g}(\xi ))=\overline{p}_{\lambda }(\xi )
\end{equation*}%
for all $g\in G$ and for all $\xi \in E$.

Let $S(G,A)=\{p\in S(A);\overline{p}_{E}\left( u_{g}\left( \xi \right)
\right) =\overline{p}_{E}\left( \xi \right) $ for all $g\in \}$. From these
facts, we conclude that $g\mapsto u_{g}$ is an inverse limit action of $G$
on $E,$ if $S(G,A)$ is a cofinal subset of $S(A).$ Therefore, if $g\mapsto
u_{g}$ is an inverse limit action of $G$ on $E$, we can suppose that $%
u_{g}=\lim\limits_{\underset{p\in S(A)}{\leftarrow }}u_{g}^{p}.$ Moreover,
the inverse limit action $g\mapsto u_{g}$ of $G$ on $E$ is continuous if and
only if the actions $g\mapsto u_{g}^{p}$ of $G$ on $E_{p},$ $p\in S(A)$ are
all continuous.

\begin{definition}
(\cite{11}) An action of $G$\ on $A$\ is a morphism $\alpha $ \ from $G$\ to
Aut$\left( A\right) $, the set of all isomorphisms of locally $C^{\ast }$%
-algebras from $A$\ to $A$. The action $\alpha $\ is continuous if the
function $\left( g,a\right) \mapsto \alpha _{g}(a)$\ from $G\times A$\ to $A$%
\ is jointly continuous.\smallskip \smallskip

A continuous action $\alpha $\ of $G$ on $A$ is an inverse limit action if
we can write $A$\ as inverse limit $\lim\limits_{\underset{\delta \in \Delta 
}{\leftarrow }}A_{\delta }$\ of $C^{\ast }$-algebras in such a way that
there are actions $\alpha ^{\delta }$\ of $G$\ on $A_{\delta }$\ such that $%
\alpha _{g}=\lim\limits_{\underset{\delta \in \Delta }{\leftarrow }}\alpha
_{g}^{\delta }$\ for all $g$\ in $G.$
\end{definition}

\begin{proposition}
(\cite{7} ) Let $G$ be a locally compact group and let $E$ be a full Hilbert 
$A$ -module. Then any action $g\mapsto u_{g}$ of $G$ on $E$ induces an
action $g\mapsto \alpha _{g}^{u}$ of $G$ on $A$ such that 
\begin{equation*}
\alpha _{g}^{u}\left( \left\langle \xi ,\eta \right\rangle \right)
=\left\langle u_{g}\left( \xi \right) ,u_{g}\left( \eta \right) \right\rangle
\end{equation*}%
for all $g\in G$ and for all $\xi ,\eta \in E$ and an action $g\mapsto \beta
_{g}^{u}$ of $G$ on $K(E)$ such that 
\begin{equation*}
\beta _{g}^{u}\left( \theta _{\xi ,\eta }\right) =\theta _{u_{g}\left( \xi
\right) ,u_{g}\left( \eta \right) }
\end{equation*}%
for all $g\in G$ and for all $\xi ,\eta \in E$. Moreover, if $g\mapsto u_{g}$
is a continuous inverse limit action of $G$ on $E,$ then the actions of $G$
on $A$ respectively $K(E)$ induced by $u$ are continuous inverse limit
actions.
\end{proposition}

\section{Action on the linking locally $C^{\ast }$-algebra of a Hilbert
module}

\begin{proposition}
Let $G$ be a locally compact group, let $E$ be a full Hilbert $A$-module.
Any action $g\mapsto u_{g}$ of $G$ on $E$ induces an action $g\mapsto \gamma
_{g}^{u}$ of $G$ on the linking algebra $\mathcal{L}(E)$ of $E$ such that 
\begin{equation*}
\gamma _{g}^{u}\left( L_{a,\xi ,\eta ,T}\right) =L_{\alpha _{g}^{u}\left(
a\right) ,u_{g}\left( \xi \right) ,u_{g}\left( \eta \right) ,\beta
_{g}^{u}\left( T\right) }
\end{equation*}%
for all $a\in A,$ $\xi ,\eta \in E$ and $T\in K(E)$. Moreover, if $g\mapsto
u_{g}$ is a continuous inverse limit action, then $g\mapsto \gamma _{g}^{u}$
is a continuous inverse limit action.
\end{proposition}

\begin{proof}
Let $g\in G$. The map $w_{g}^{u}:A\oplus E\rightarrow A\oplus E$ defined by 
\begin{equation*}
w_{g}^{u}\left( a\oplus \xi \right) =\alpha _{g}^{u}(a)\oplus u_{g}\left(
\xi \right)
\end{equation*}%
is a morphism of Hilbert modules, since 
\begin{eqnarray*}
\left\langle w_{g}^{u}\left( a\oplus \xi \right) ,w_{g}^{u}\left( b\oplus
\eta \right) \right\rangle &=&\left\langle \alpha _{g}^{u}(a),\alpha
_{g}^{u}(b)\right\rangle +\left\langle u_{g}\left( \xi \right) ,u_{g}\left(
\eta \right) \right\rangle \\
&=&\alpha _{g}^{u}\left( \left\langle a\oplus \xi ,b\oplus \eta
\right\rangle \right)
\end{eqnarray*}%
for all $a,b\in A$ and for all $\xi ,\eta \in E$, and $\alpha _{g}^{u}$ is
an isomorphism of locally $C^{\ast }$-algebras. Moreover, since $w_{g}^{u}$
is invertible and $\left( w_{g}^{u}\right) ^{-1}=w_{g^{-1}}^{u}$, $w_{g}^{u}$
is an isomorphism of Hilbert modules. It is not difficult to check that $%
g\mapsto w_{g}^{u}$ is an action of $G$ on $A\oplus E.$

Let $\gamma ^{u}$ be the action of $G$ on $K(A\oplus E)$ induced by $w^{u}$.
Then 
\begin{equation*}
\gamma _{g}^{u}\left( \theta _{a\oplus \xi ,b\oplus \eta }\right) =\theta
_{w_{g}^{u}\left( a\oplus \xi \right) ,w_{g}^{u}\left( b\oplus \eta \right)
}=\theta _{\alpha _{g}^{u}(a)\oplus u_{g}\left( \xi \right) ,\alpha
_{g}^{u}(b)\oplus u_{g}\left( \eta \right) }
\end{equation*}%
for all $a,b\in A$, for all $\xi ,\eta \in E$, and for all $g\in G$.

Let $g\in G,$ $a\in A,$ $\xi ,\eta \in E$ and $T\in K(E).$ We will show that 
\begin{equation*}
\gamma _{g}^{u}\left( L_{a,\xi ,\eta ,T}\right) =L_{\alpha _{g}^{u}\left(
a\right) ,u_{g}\left( \xi \right) ,u_{g}\left( \eta \right) ,\beta
_{g}^{u}\left( T\right) }.
\end{equation*}%
For this, let $\{e_{i}\}_{i}$ be an approximate unit for $A.$ From 
\begin{equation*}
\widetilde{p}_{L(A\oplus E)}(L_{a,0,0,0}-\theta _{a\oplus 0,e_{i}\oplus
0})\leq p\left( a-ae_{i}\right)
\end{equation*}%
\begin{equation*}
\widetilde{p}_{L(A\oplus E)}\left( L_{0,\xi ,0,0}-\theta _{e_{i}\oplus
0,0\oplus \xi }\right) \leq \overline{p}_{E}\left( \xi -\xi e_{i}\right)
\end{equation*}%
and 
\begin{equation*}
\widetilde{p}_{L(A\oplus E)}\left( L_{0,0,\eta ,0}-\theta _{0\oplus \eta
,e_{i}\oplus 0}\right) \leq \overline{p}_{E}\left( \eta -\eta e_{i}\right)
\end{equation*}
for all $p\in S(A)$ and for all $i\in I$, and taking into account that $%
\gamma _{g}^{u},$ $\alpha _{g}^{u}$ and $u_{g}$ are continuous, $%
ae_{i}\rightarrow a,$ $\xi e_{i}\rightarrow \xi $ and $\eta e_{i}\rightarrow
\eta ,$ we conclude that 
\begin{equation*}
\gamma _{g}^{u}\left( L_{a,\xi ,\eta ,0}\right) =L_{\alpha _{g}^{u}\left(
a\right) ,u_{g}\left( \xi \right) ,u_{g}\left( \eta \right) ,0}\text{.}
\end{equation*}%
If $T\in K(E)$, then there is a net $\{\sum\limits_{k\in I_{j}}\theta _{\xi
_{k},\eta _{k}}\}_{j}$ in $\Theta (E)$ which converges to $T$. From 
\begin{equation*}
\widetilde{p}_{L(A\oplus E)}(L_{0,0,0,T}-\sum\limits_{k\in I_{j}}\theta
_{0\oplus \xi _{k},0\oplus \eta _{k}})\leq \widetilde{p}_{L(E)}(T-\sum%
\limits_{k\in I_{j}}\theta _{\xi _{k},\eta _{k}})
\end{equation*}%
for all $p\in S(A)$, and taking into account that $\gamma _{g}^{u}$ and $%
\beta _{g}^{u}$ are continuous and 
\begin{equation*}
\sum\limits_{k\in I_{j}}\gamma _{g}^{u}\left( \theta _{0\oplus \xi
_{k},0\oplus \eta _{k}}\right) =\sum\limits_{k\in I_{j}}\theta _{0\oplus
u_{g}(\xi _{k}),0\oplus u_{g}(\eta _{k})}=0\oplus \beta
_{g}^{u}(\sum\limits_{k\in I_{j}}\theta _{\xi _{j},\eta _{j}})
\end{equation*}%
we deduce that 
\begin{equation*}
\gamma _{g}^{u}\left( L_{0,0,0,T}\right) =L_{0,0,0,\beta _{g}^{u}\left(
T\right) }.
\end{equation*}%
Thus we have 
\begin{eqnarray*}
\gamma _{g}^{u}\left( L_{a,\xi ,\eta ,T}\right) &=&\gamma _{g}^{u}\left(
L_{a,\xi ,\eta ,0}\right) +\gamma _{g}^{u}\left( L_{0,0,0,T}\right) \\
&=&L_{\alpha _{g}^{u}\left( a\right) ,u_{g}\left( \xi \right) ,u_{g}\left(
\eta \right) ,0}+L_{0,0,0,\beta _{g}^{u}\left( T\right) } \\
&=&L_{\alpha _{g}^{u}\left( a\right) ,u_{g}\left( \xi \right) ,u_{g}\left(
\eta \right) ,\beta _{g}^{u}\left( T\right) }.
\end{eqnarray*}

If $u$ is a continuous inverse limit action, then we can suppose that $%
u_{g}=\lim\limits_{\underset{p\in S(A)}{\leftarrow }}u_{g}^{p},$ where $%
g\mapsto u_{g}^{p}$ is a continuous action of $G$ on $E_{p}$ for each $p\in
S(A).$ Let $p\in S(A)$ and let $g\mapsto w_{g}^{u^{p}}$ be the action of $G$
on $A_{p}\oplus E_{p}$ induced by $u^{p}$. It is not difficult to check that 
$\left( w_{g}^{u^{p}}\right) _{p}$ is an inverse system of isomorphisms of
Hilbert $C^{\ast }$-modules and $g\mapsto \lim\limits_{\underset{p}{%
\leftarrow }}w_{g}^{u^{p}}$ is a continuous inverse limit action of $G$ on $%
A\oplus E.$ Moreover, $w_{g}^{u}=\lim\limits_{\underset{p}{\leftarrow }%
}w_{g}^{u^{p}}$ for each $g\in G$. By Proposition 2.5, the action $\gamma
^{u}$ of $G$ on $K(A\oplus E)$ induced by $w^{u}$ is a continuous inverse
limit action, and moreover, $\gamma _{g}^{u}=\lim\limits_{\underset{p}{%
\leftarrow }}\gamma _{g}^{u^{p}}$ for each $g\in G$, where $\gamma ^{u^{p}}$
is the action of $G$ of $\mathcal{L}(E_{p})$ induced by $u^{p}$.
\end{proof}

\begin{remark}
Let $G$ be a locally compact group, let $E$ be a full Hilbert $A$-module,
and let $g\mapsto u_{g}$ be an action of $G$ on $E$.

\begin{enumerate}
\item Since the map $a\mapsto L_{a,0,0,0}$ from $A$ to $\mathcal{L}(E)$
identifies $A$ with a locally $C^{\ast }$-subalgebra of $\mathcal{L}(E)$ and 
\begin{equation*}
\gamma _{g}^{u}\left( L_{a,0,0,0}\right) =L_{\alpha _{g}^{u}\left( a\right)
,0,0,0}
\end{equation*}%
for all $a\in A$ and for all $g\in G$, the restriction of $\gamma ^{u}$ to $%
A $ can be identified with the action of $G$ on $A$ induced by $u$.

\item Since the map $T\mapsto L_{0,0,0,T}$ from $K(E)$ to $\mathcal{L}(E)$
identifies $K(E)$ with a locally $C^{\ast }$-subalgebra of $\mathcal{L}(E)$
and 
\begin{equation*}
\gamma _{g}^{u}\left( L_{0,0,0,T}\right) =L_{0,0,0,\beta _{g}^{u}(T)}
\end{equation*}%
for all $T\in K(E)$ and for all $g\in G$, the restriction of $\gamma ^{u}$
to $K(E)$ can be identified with the action of $G$ on $K(E)$ induced by $u$.
\end{enumerate}
\end{remark}

\begin{definition}
(\cite{7}) Let $G$ be a locally compact group and let $g\mapsto \alpha _{g}$
and $g\mapsto \beta _{g}$ be two continuous inverse limit actions of $G$ on
two locally $C^{\ast }$-algebras $A$ and $B$. We say that the actions $%
\alpha $ and $\beta $ are conjugate if there is a an isomorphism of locally $%
C^{\ast }$-algebras $\varphi :A\rightarrow B$ such that $\alpha _{g}=\varphi
^{-1}\circ \beta _{g}\circ $ $\varphi $ for each $g\in G.$
\end{definition}

\begin{definition}
(\cite{7}) Let $G$ be a locally compact group and let $g\mapsto \alpha _{g}$
and $g\mapsto \beta _{g}$ be two continuous inverse limit actions of $G$ on
two locally $C^{\ast }$-algebras $A$ and $B$. We say that $\alpha $ and $%
\beta $ are strongly Morita equivalent if there is a full Hilbert module $E$
over $A$, and there is a continuous inverse limit action $g\mapsto u_{g}$ of 
$G$ on $E$ such that the actions of $G$ on $A$ respectively $K(E)$ induced
by $u$ are conjugate with $\alpha $ respectively $\beta .$
\end{definition}

Recall that two locally $C^{\ast }$-algebras $A$ and $B$ are two
complementary corners in a given locally $C^{\ast }$-algebra $C$, if there
is two projections $e$ and $f$ in the multiplier algebra $M(C)$ of $C$ such
that:

\begin{enumerate}
\item $A=eCe$ and $B=fCf$;

\item $e+f=1_{M(C)};$

\item the locally $C^{\ast }$-subalgebras $CeC$ and $CfC$ of $C$ are dense
in $C$.
\end{enumerate}

The following theorem is a version of Theorem 2.9 \cite{5}.

\begin{theorem}
Let $G$ be a locally compact group and let $g\mapsto \alpha _{g}$ and $%
g\mapsto \beta _{g}$ be two continuous inverse limit actions of $G$ on two
locally $C^{\ast }$-algebras $A$ and $B$. Then the actions $\alpha $ and $%
\beta $ are strongly Morita equivalent if and only if there is a locally $%
C^{\ast }$-algebra $C$ such that $A$ and $B$ appear as two complementary
full corners in $C$ and there is a continuous inverse limit action $g\mapsto
\gamma _{g}$ of $G$ on $C$ such that $A$ and $B$ are invariant to $\gamma $
and the actions $g\mapsto \gamma _{g}|_{A}$ and $g\mapsto \gamma _{g}|_{B}$
of $G$ on $A$ respectively $B$ can be identified with $\alpha $ respectively 
$\beta .$
\end{theorem}

\begin{proof}
First we suppose that $\alpha $ and $\beta $ are strongly Morita equivalent.
Let $(E,u)$ be the pair consisting of a full Hilbert $A$-module and a
continuous inverse limit action of $G$ on $E$ which implements a strong
Morita equivalence between $\alpha $ and $\beta $. Let $C=\mathcal{L}(E)$
and let $\gamma ^{u}$ be the action of $G$ on $C$ induced by $u$. Then $A$
and $B$ are isomorphic with two complementary full corners in $C$ ( Theorem
2.9 in \cite{5}) and $g\mapsto \gamma _{g}^{u}$ is a continuous inverse
limit action of $G$ on $C$ such that identifying $A$ and $B$ with corners in 
$C$, $\gamma ^{u}|_{A}=\alpha $ and $\gamma ^{u}|_{B}=\beta $ (Remark 3.2).

Conversely, we suppose that there is a locally $C^{\ast }$-algebra $C$ such
that $A$ and $B$ appear as two complementary full corners in $C$ and there
is a continuous inverse limit action $g\mapsto \gamma _{g}$ of $G$ on $C$
such that $A$ and $B$ are invariant to $\gamma $ and the actions $g\mapsto
\gamma _{g}|_{A}$ and $g\mapsto \gamma _{g}|_{B}$ of $G$ on $A$ respectively 
$B$ can be identified with $\alpha $ respectively $\beta $. By Proposition
2.8 in \cite{5}, the locally $C^{\ast }$-algebras $A$ and $C$ are strongly
Morita equivalent. Moreover, if $e$ is a full projection in $M(C)$, the
multiplier algebra of $C$, such that $A=eCe$, then the Hilbert $A$-module $%
Ce $ implements a strong Morita equivalence between $A$ and $C$. Let $g\in G$%
. For each $c\in C,$ $\gamma _{g}(ce)\in Ce$, since $\gamma _{g}(ce)e=\gamma
_{g}(ce)$. Thus we can consider the linear map $u_{g}:Ce\rightarrow Ce$
defined by $u_{g}(ce)=\gamma _{g}(ce)$. Since 
\begin{eqnarray*}
\left\langle u_{g}(ce),u_{g}(de)\right\rangle &=&\left\langle \gamma
_{g}(ce),\gamma _{g}(de)\right\rangle =\gamma _{g}\left( ec^{\ast }de\right)
\\
&=&\alpha _{g}\left( ec^{\ast }de\right) =\alpha _{g}\left( \left\langle
ce,de\right\rangle \right)
\end{eqnarray*}%
for all $c,d\in C$ and since $u_{g}$ is invertible and $\left( u_{g}\right)
^{-1}=u_{g^{-1}},$ $u_{g}\in $Aut$(Ce)$. It is not difficult to check that $%
g\mapsto u_{g}$ is an action of $G$ on $Ce.$ Moreover, since $\gamma $ is a
continuous inverse limit action, $u$ is a continuous inverse limit action.

Since 
\begin{equation*}
\left\langle u_{g}(ce),u_{g}(de)\right\rangle =\alpha _{g}\left(
\left\langle ce,de\right\rangle \right)
\end{equation*}%
for all $g\in G$ and for all $c,d\in C,$ $\alpha ^{u}=\alpha $. From 
\begin{equation*}
\beta _{g}^{u}(\theta _{ce,d^{\ast }e})=\theta _{u_{g}(ce),u_{g}(d^{\ast
}e)}=\theta _{\gamma _{g}(ce),\gamma _{g}(d^{\ast }e)}
\end{equation*}%
for all $g\in G$ and for all $c,d\in C,$ and taking into account that $CeC$
is dense in $C$ and an element $ced$ in $C$ can be identified with the
element $\theta _{ce,d^{\ast }e}$ in $K(Ce),$ we deduce that the actions $%
\beta ^{u}$ and $\gamma $ are conjugate. Thus we proved that the actions $%
\alpha $ and $\gamma $ are strongly Morita equivalent. In the same way we
show that the actions $\beta $ and $\gamma $ are equivalent, and since the
strong Morita equivalence is an equivalence relation \cite{7}, the actions $%
\alpha $ and $\beta $ are strongly Morita equivalent.
\end{proof}

Using Lemma 5.2 in \cite{11} and Theorem 3.5 we obtain the following
corollary.

\begin{corollary}
Let $G$ be a compact group and let $g\mapsto \alpha _{g}$ and $g\mapsto
\beta _{g}$ be two continuous actions of $G$ on two locally $C^{\ast }$%
-algebras $A$ and $B.$ Then the actions $\alpha $ and $\beta $ are strongly
Morita equivalent if and only if there is a locally $C^{\ast }$-algebra $C$
such that $A$ and $B$ appear as two complementary full corners in $C$ and
there is a continuous inverse limit action $g\mapsto \gamma _{g}$ of $G$ on $%
C$ such that $A$ and $B$ are invariant to $\gamma $ and the actions $%
g\mapsto \gamma _{g}|_{A}$ and $g\mapsto \gamma _{g}|_{B}$ of $G$ on $A$
respectively $B$ can be identified with $\alpha $ respectively $\beta .$
\end{corollary}

Department of Mathematics, Faculty of Chemistry, University of Bucharest,
Bd. Regina Elisabeta nr. 4-12, Bucharest, Romania

mjoita@fmi.unibuc.ro

\bigskip

\end{document}